\documentclass[aop,MSNbibl,dvips]{arximspdf}
\usepackage{graphicx}

\doi{10.1214/12-AOP784} %kopijuoti is PTS
\volume{41}
\issue{3B}
\pubyear{2013}
\firstpage{2066}
\lastpage{2089}

\makeatletter

\newcommand{\ovverline}{\bar}

\newcommand{\rrvert}{\vert}
\newcommand{\llvert}{\vert}

\newcommand{\xrightarrow}[1]{\mathop{\hbox to 1cm{\rightarrowfill}}_{#1}}

\newtheorem{theorem}{Theorem}[section]
\newtheorem{proposition}[theorem]{Proposition}
\newtheorem{lemma}[theorem]{Lemma}

\newproclaim{definition}[theorem]{Definition}

\newcommand{\bbR}{\mathbb{R}}
\newcommand{\E}{\mathbb{E}}
\newcommand{\R}{\bbR}
\newcommand{\bg}{\bar{g}}

\newcommand{\bbP}{\mathbb{P}}

\newcommand{\const}{\operatorname{const}}
\newcommand{\supp}{\operatorname{supp}}
\newcommand{\mes}{\operatorname{mes}}

\newcommand{\card}{\operatorname{card}}

\makeatother

\begin{document}
\begin{frontmatter}

\title{Symmetric random walks on $\mathbf{Homeo}^+(\mathbf{R})$}
\runtitle{Symmetric random walks on $\mathrm{Homeo}^+(\mathbf{R})$}

\begin{aug}
\author[A]{\fnms{B.} \snm{Deroin}\thanksref{t1}\ead[label=e1]{Bertrand.Deroin@math.u-psud.fr}},
\author[B]{\fnms{V.} \snm{Kleptsyn}\corref{}\thanksref{t2}\ead[label=e2]{Victor.Kleptsyn@univ-rennes1.fr}},
\author[C]{\fnms{A.} \snm{Navas}\thanksref{t3}\ead[label=e3]{andres.navas@usach.cl}}
\and
\author[D]{\fnms{K.} \snm{Parwani}\thanksref{t4}\ead[label=e4]{kparwani@eiu.edu}}
\runauthor{Deroin, Kleptsyn, Navas and Parwani}
\affiliation{Universit\'e Paris-Sud, Institut de Recherches Math\'
ematiques de Rennes,
Universidad de Santiago de Chile and Eastern Illinois University}
\dedicated{Dedicated to John Milnor on his 80th anniversary}
\address[A]{B. Deroin\\
Lab. de Math\'ematiques\\
CNRS, Universit\'e Paris-Sud\\
B\^at 425\\
91405 Orsay Cedex\\
France\\
\printead{e1}} %adresu isvedimo komanda gale!
\address[B]{V. Kleptsyn\\
Institut de Recherches Math\'ematiques\\
\quad de Rennes (UMR 6625 CNRS)\\
Campus Beaulieu, 35042 Rennes\\
France\\
\printead{e2}}
\address[C]{A. Navas\\
Universidad de Santiago de Chile\\
Alameda 3363, Santiago\\
Chile\\
\printead{e3}}
\address[D]{K. Parwani\\
Department of Mathematics\hspace*{44pt}\\
Eastern Illinois University\\
600 Lincoln Avenue\\
Charleston, Illinois 61920\\
USA\\
\printead{e4}}
\end{aug}

\thankstext{t1}{Supported by ANR-08-JCJC-0130-01 and ANR-09-BLAN-0116.}

\thankstext{t2}{Supported by ANR Symplexe BLAN-06-3-137-237, the RFBR
Grant 10-01-00739-a
and the RFBR/CNRS Grant 10-01-93115-CNRS-a.}

\thankstext{t3}{Supported by the Conicyt-PBCT Grant ADI 17 and the
Fondecyt Grant
1100536.}

\thankstext{t4}{Supported in part by ANR symplexe BLAN-06-3-137-237
(France).}

% HISTORY:
\received{\smonth{3} \syear{2011}}
\revised{\smonth{3} \syear{2012}}

% ABSTRACT
%
\begin{abstract}
We study symmetric random walks on finitely generated
groups of orientation-preserving homeomorphisms of the real line. We establish
an oscillation property for the induced Markov chain on the line that
implies a weak form of recurrence. Except for a few special cases, which
can be treated separately, we prove a property of ``global stability at a
finite distance'': roughly speaking, there exists a compact interval such
that any two trajectories get closer and closer whenever one of them
returns to the compact interval. The probabilistic techniques employed
here lead to interesting results for the study of group actions on the
line. For instance, we show that under a suitable change of the coordinates,
the drift of every point becomes zero provided that the action is
minimal. As
a byproduct, we recover the fact that every finitely generated group of
homeomorphisms
of the real line is topologically conjugate to a group of (globally)
Lipschitz homeomorphisms.
Moreover, we show that such a conjugacy may be chosen in such a way
that the displacement of
each element is uniformly bounded.
\end{abstract}

% KEYWORDS

\end{frontmatter}

%s1 #&#
\section{Introduction}

In this article, we study symmetric random walks on
finitely generated groups
of (orientation-preserving) homeomorphisms of the real line. The
results presented here fit into the general framework
of systems of iterated random functions~\cite{IRF}. However, besides
the lack of compactness of
the phase space, there is a crucial point that separates our approach
from the classical
ones---the complete absence of any hypothesis of contraction. So to
carry out our study, we
need to use some extra structure, and this is provided by the natural
ordering of the
real line. In this direction, the results herein are also closely
related to~\cite{kellerer}, where general
Markov processes on ordered spaces are examined. However, since we only
consider symmetric
measures, there is no zero drift condition required for our processes,
unlike~\cite{kellerer} where this
is a crucial assumption. In fact, one of our main results
is that when the action is minimal, an appropriate change of
coordinates on the
real line makes the drift of \textit{every point} equal to zero. This
follows from a
reparametrization that utilizes the stationary measure and a
straightforward argument that
employs the one-dimensional structure of the phase space in a decisive manner.

It is quite remarkable that, in presence of the linear order structure,
we recover several
phenomena that in more complex phase spaces are very specific to
particular classes of groups.
For instance, in~\cite{BBE}, Babillot, Bougerol and Elie consider
random walks
on the group of affine homeomorphisms of ${\mathbf R}^n$ in the difficult
case where the
logarithm of the expansion rate vanishes in mean. In this situation,
they show the existence
of an infinite Radon measure that is invariant by the transition
operator, and for the
case where the Lebesgue measure is not totally invariant, they
establish a property of
``global stability at a finite distance'': any two trajectories get
closer and closer whenever
one of them returns to a fixed compact set. Using this property, they
obtain the uniqueness
of the stationary measure (up to a constant multiple). It turns out
that three of
our main results here are analogues of these facts for groups
of homeomorphisms of the real line. Furthermore,
these results are also analogous to---though much more
elaborate than---the previously established results for groups of
circle homeomorphisms; see, for instance,~\cite{DKN}, Section 5.1.
As in the case of~\cite{DKN}, the proofs here involve a prior study of the
general structure of the associated dynamics, which is the core of this paper.

The motivation for studying groups of homeomorphisms of the real line
comes from many
sources. Algebraically, these groups are characterized by the existence
of a left-invariant
total order~\cite{book}, which fits into well developed and quite
formal theories~\cite{glass}.
More recently, many results about groups acting on the real line or the
circle have focused on the
relation with ``rigidity theory,'' a kind of nonlinear version of
representation theory where one seeks
to understand the nature of the obstructions to the existence of
(faithful) group
actions on specific phase spaces (see~\cite{fisher} for a survey of
these ideas). In this direction, it is
conjectured that some particular groups, like groups with Kazhdan's
property (T) or lattices in higher-rank simple Lie groups, do not act
on the real line (or equivalently,
are not left-orderable). We strongly believe that our probabilistic
approach opens
new and promising avenues of study that bear the potential to yield
important results
in the investigation of these and many other open questions concerning
left-orderable groups.

%%%%%%%%%%%%%%%%%%%%%%%%%%%%%%%%%%%%%%%%%%%%%%%%%%%%%%%%%%%%%%%%%%%%%%%%%%%%%%%%%%%%%%%%%%%%%%%%%%

%s2 #&#
\section{A description of the results}

Given a symmetric probability measure~$\mu$, we
consider the group $G$ generated by its support. Although some of our
results apply to the case where this support is countably infinite,
we restrict our discussion to the case where it is finite.
Assume throughout that the action
of $G$ is \textit{irreducible}; that is, there are no global fixed
points. Otherwise,
one may consider the action the connected components of the complement
of the set
of these global fixed points (on each of these components,
the action is irreducible). We then consider the random walk
induced by $\mu$ as a Markov process on the real line. In Section \ref
{Srecurrence},
the recurrence properties of this process are studied. We first prove
that almost every
trajectory oscillates between $-\infty$ and $+\infty$ (Proposition
\ref
{Poscillation}).
Moreover, there exists a compact interval $K$ such that these
trajectories pass
through $K$ infinitely many times (Theorem~\ref{Crecurrence}). Using standard
arguments a la Chacon--Orstein, this allows us to show the existence of
a stationary
Radon measure on the line (Theorem~\ref{Cexistenceofinvariantmeasure}).

In Section~\ref{Spropertiesofinvariantmeasures}, general properties
of the
stationary measures are examined. If the atomic part of the stationary
measure is nontrivial, then it is supported on the union of discrete
orbits and is totally invariant (Lemmas~\ref{Latomicpart1} and~\ref{Latomicpart2}).
If there are no discrete orbits, there exists
a unique minimal nonempty closed invariant set $\mathcal{M}$ that is
the support of any stationary measure (Proposition~\ref{Pminimal} and
Lemma~\ref{Lsupport} ). Furthermore,
the stationary measure is unique up to a constant factor (Theorem
\ref{Tunicity}). This result follows from an argument
due to S. Brofferio in~\cite{Brofferio} and a nondivergence property
for the
trajectories of the Markov process established in Lemma~\ref{Lweakcontraction}.

In Section~\ref{Sstability}, we obtain the property of global stability at
a finite distance
provided that no invariant Radon measure exists and $G$ is not
centralized by any
homeomorphism without fixed points. Roughly speaking, this last
condition means that the action
does not appear as the lift of an action on the circle. If this is not
the case, a weak
form of the contraction property is established (all of this is
summarized in
Theorem~\ref{Tprobabilisticcontractionproperty}).

In Section~\ref{Sderriennic}, we provide a connection to the
beautiful work~\cite{Derriennic},
where Derriennic studies Markov processes on the real line satisfying
$\mathbb E( X^x_1 ) = x$
for large values of $|x|$.
For every finitely generated group of homeomorphisms of the real line
acting minimally,
we produce a coordinate change for which the Derriennic property
[$\mathbb E( X^x_1 ) = x$] holds
for \textit{every} $x \in{\mathbf R}$ (Theorem \ref
{Tderrienniccoordinates}). This is done by appropriately
integrating the associated stationary measure. A careful analysis of
the invariant Radon measure is
carried out before establishing this result. In particular, we prove
that the measure is infinite
on every unbounded interval (Lemma~\ref{Lbi-infiniteness}). As a
consequence of the existence
of these Derriennic coordinates, we recover a rather surprising fact:
every finitely
generated group of homeomorphisms of the real line is topologically
conjugate to a
group of Lipschitz homeomorphisms (Theorem~\ref{lip}). (This result
also follows from
the (probabilistic) techniques introduced in~\cite{DKN}.) Moreover, we
show that such
a conjugacy may be taken so that the displacement function $x \mapsto
g(x) - x$
becomes bounded uniformly in $x$ for all $g \in G$.

%%%%%%%%%%%%%%%%%%%%%%%%%%%%%%%%%%%%%%%%%%%%%%%%%%%%%%%%%%%%%%%%%%%%%%%%%%%%%%%%%%%%%%%%%%%%%%%%%%%%%%%%%

%s3 #&#
\section{Notation}

Let $\{g_n\}$ be a sequence
of i.i.d. $\mathrm{Homeo}^+({\mathbf R})$-valued random variables,
whose distribution is a \textit{symmetric} measure $\mu$.
The left random walk on $\mathrm{Homeo}^+ ({\mathbf R})$
is defined by the random variables
\[
f_n:= g_n \circ\cdots\circ g_1.
\]
More precisely,\vspace*{1pt} let $G$ be the group generated by the support of $\mu$
and consider the probability
space $\Omega:= (G^{\mathbf N}, \mu^{\otimes{\mathbf N}} )$. Then $g_n$ is
defined to be the
$n$th coordinate of $\omega\in\Omega$. The group
$G$ is assumed to be countable, and in fact, $G$ will be finitely
generated in most cases.

We introduce the Markov chain $X$ on the real line, that is,
for any $x\in{\mathbf R}$ and any nonnegative integer $n$,
\[
X_n^x:= f_n (x).
\]
Let $C_b({\mathbf R})$ and $C_c({\mathbf R})$ denote the spaces of
continuous bounded functions and compactly supported continuous
functions, respectively. Let
$P\dvtx C_b({\mathbf R})\rightarrow C_b({\mathbf R})$ be the
transition operator,
defined as usual by
\[
P (\varphi) (x):= \mathbb E \bigl(\varphi\bigl(X_1^x
\bigr)\bigr),
\]
where $\varphi\in C_b({\mathbf R})$ and $x \in{\mathbf R}$. The
operator $P$
acts by duality on the set of finite measures on the real line,
and if $\mu$ is finitely supported, $P$ preserves $C_c({\mathbf R})$
and hence acts by the duality on the set of Radon measures.

It should come as no surprise that in the investigation of a group
action on the real line, one
is led to study the action on the components of ${\mathbf R} \setminus
\operatorname{Fix}(G)$.
On any such component, no global fixed point exists, and so, we may
(and we always will) assume that the
action is \textit{irreducible}, that is, for every $x \in{\mathbf R}$ there
exists a
$g \in G$ such that $g(x) \neq x$. If $\mu$ is a symmetric probability
measure on
$\mathrm{Homeo}^+ ({\mathbf R})$, then we will say that $\mu$ is irreducible
if the group generated by its support satisfies this property.

%%%%%%%%%%%%%%%%%%%%%%%%%%%%%%%%%%%%%%%%%%%%%%%%%%%%%%%%%%%%%%%%%%%%%%%%%%%%%%%%%%%%%%%%%%%%%%%%%%%

%s4 #&#
\section{Recurrence} \label{Srecurrence}

In this section, we establish the recurrence of the
Markov chain $X$
when the measure $\mu$ is irreducible, symmetric, and has finite
support. We begin
with a lemma that extends~\cite{DKN}, Proposition 5.7 (see also \cite
{FFP}) and that
will be crucial in Sections~\ref{Suniqueness} and~\ref{Sderriennic}.
The proof
is based on the second proof proposed in~\cite{DKN}; for a proof based
on the first proof therein, see~\cite{book}.
%
%le4.1 #&#
%
\begin{lemma} \label{Lbi-infiniteness}
Let $\mu$ be an irreducible, symmetric probability measure
on the group $\mathrm{Homeo}^+ ({\mathbf R})$. Then any nonvanishing
$P$-invariant
Radon measure $\nu$ on the real line is bi-infinite [i.e., $\nu
(x,\infty
) = \infty$
and $\nu(-\infty,x) =\infty$, for all $x \in{\mathbf R}$].
\end{lemma}
\begin{pf}
Suppose that there exists an $x\in{\mathbf R}$ such that $\nu(x,\infty) <
\infty$.
Since the action is irreducible, for every $y\in{\mathbf R}$, there is an
element $g\in G$ such that $g(x) < y$.
Select $n > 0$ such that $\mu^{\star n} (g^{-1}) >0$ and then observe that
\[
\nu(y,\infty) \leq\nu\bigl(g(x),\infty\bigr) \leq\frac{1}{\mu^{\star
n}(g^{-1})} \nu(x,
\infty)<\infty.
\]
This argument implies that $\nu(y,\infty) <\infty$ for all $y \in
{\mathbf R}$.

Now let $f\dvtx{\mathbf R} \rightarrow(0,\infty)$ be the function
defined by
$f(x):= \nu(x,\infty)$.
Since $\mu$ is symmetric, this function is harmonic on the orbits.
Fix a real number $A$ satisfying $0<A< \nu(-\infty,\infty)$, and then
let $h:= \max(0, A - f)$. The function $h$
is subharmonic, that is, $h \leq P(h)$. Moreover, it vanishes on a
neighborhood of $-\infty$
and is bounded on a neighborhood of $\infty$.
This implies that $h$ is $\nu$-integrable, and since $\int Ph \,d\nu=
\int h \,d\nu$,
the function $h$ must be $P$-invariant $\nu$-a.e. Now a classical lemma
in~\cite{Garnett} asserts that a measurable
function which is in $L^1(\mathbf R,\nu)$ and $P$-invariant must be
$G$-invariant almost everywhere.
Thus, $h$ is constant on almost every orbit. However, this is
impossible since every
orbit intersects every neighborhood of $-\infty$ (where $h$ vanishes)
and of $\infty$ (where $h$ is positive). This contradiction establishes
the desired result.
\end{pf}

%pr4.2 #&#
%
\begin{proposition}[(Oscillation)]
\label{Poscillation}
Let $\mu$ be an irreducible, symmetric probability measure on $\mathrm
{Homeo}^+ ({\mathbf R})$.
Then for every $x \in{\mathbf R}$, almost surely we have
\[
\limsup_{n\rightarrow\infty} X_n^x = +\infty\quad\mbox{and}\quad
\liminf_{n\rightarrow\infty} X_n^x = -\infty.
\]
\end{proposition}
\begin{pf}
Given points $A$ and $x$ on the real line, let
\[
p_A (x):= \mathbb P \Bigl\{ \limsup_{n\rightarrow\infty}
X_n^x > A \Bigr\}.
\]
Since $G$ acts by orientation-preserving homeomorphisms,
for each $x \leq y$ we have
\[
\Bigl\{(g_n)_n\in G^{\mathbf N} \bigm|
\limsup_{n\rightarrow\infty} X_n^x > A \Bigr\} \subset\Bigl
\{(g_n)_n\in G^{\mathbf N} \bigm|\limsup_{n\rightarrow\infty}
X_n^y > A \Bigr\}.
\]
In particular, $p_A(x) \leq p_A(y)$, that is, $p_A$ is nondecreasing.
Moreover, since $p_A$ is the probability of the tail event
\[
\Bigl\{\limsup_{n\rightarrow\infty} X_n^x > A \Bigr\},
\]
and $X$ is a Markov chain, $p_A$ is \textit{harmonic}:
for every $x\in{\mathbf R}$ and every integer $n\geq0$,
\[
p_A(x) = \mathbb E \bigl( p_A \bigl(X_n^x
\bigr)\bigr).
\]

We would like to think of $p_A$ as the distribution function of a
\textit{finite}
measure on ${\mathbf R }$. Since this is possible only\vadjust{\goodbreak} if $p_A$ is continuous
on the right, we are led to consider the right-continuous function
\[
\ovverline{p}_A(x): = \lim_{y\rightarrow x, y>x} p_A (y).
\]
This function is still nondecreasing; hence there exists a finite
measure $\nu$ on ${\mathbf R}$ such that for all $x < y$,
\[
\nu(x,y] = \ovverline{p}_A (y) -
\ovverline{p}_A (x).
\]
Since $p_A$ is harmonic and $G$ acts by homeomorphisms, the function
$\ovverline{p}_A$
must be harmonic. Thus $\nu$ is also harmonic, and furthermore, the
measure $\nu$ is $P$-invariant since $\mu$ is symmetric. Now recall
that Lemma~\ref{Lbi-infiniteness} implies that any $P$-invariant
finite measure vanishes
identically (see also~\cite{DKN}, Proposition~5.7), and therefore,
$\nu
=0$ and $\ovverline{p}_A$
is constant. The 0--1 law can be applied here to conclude that (for any
fixed $A$) either $p_A(x) \equiv0$ or $p_A(x) \equiv1$.

Let us now show that $p_A$ equals to $1$ for any $A$. Indeed, fix any
$x_0>A$. Then,
for any $g\in\mathrm{Homeo}^+({\mathbf R})$, we have either $g(x_0)\ge
x_0$, or $g^{-1}(x_0)\ge x_0$,
and hence, due to the symmetry of measure $\mu$, for every $n$ the
inequality $X^{x_0}_n\ge x_0$ holds with probability at least $1/2$. It
is easy then to see that
\[
p_A = p_A(x_0) \ge\limsup_{n \to\infty}
\mathbb{P} \bigl\{ X_n^{x_0} \ge x_0\bigr\}
\ge1/2.
\]
As we have already shown that $p_A$ is equal to $0$ or to $1$, this
implies that $p_A$ is identically equal to $1$.

The latter means that for every $x\in{\mathbf R}$,
\[
\limsup_{n\rightarrow\infty} X_n^x = +\infty
\]
holds almost surely.
Analogously, for every $x \in{\mathbf R}$, almost surely we have
\[
\liminf_{n\rightarrow\infty} X_n^x = -\infty.
\]
This completes the proof of the proposition.
\end{pf}

We are now ready to prove the main result of this section.
%
%th4.3 #&#
%
\begin{theorem}[(Recurrence)]
\label{Crecurrence}
Let $\mu$ be an irreducible, finitely supported, symmetric probability
measure on
$\mathrm{Homeo}^+({\mathbf R})$. Then there exists a compact interval $K$
such that, for
every $x$, almost surely the sequence $(X_n^x)$ intersects $K$
infinitely often.
\end{theorem}
\begin{pf} Consider an interval $K = [A, B]$, where $A<B$ are points in the
real line such that for every element $g$ of the support of $\mu$, we
have $g(A) < B$.
By Proposition~\ref{Poscillation}, for every $x\in{\mathbf R}$, almost
surely the
sequence $(X_n^x)$ will pass from $(-\infty, A]$ to $[B,+\infty)$
infinitely often. Now the desired conclusion follows from the
observation that the choices of $A$ and $B$ imply that every time this
happens, $\{ X_n^x\}$ must traverse the interval $K$.
\end{pf}

%%%%%%%%%%%%%%%%%%%%%%%%%%%%%%%%%%%%%%%%%%%%%%%%%%%%%%%%%%%%%%%%%%%%%%%%%%%%%%%%%%%%%%%%%%%%%%%%

%s5 #&#
\section{Existence of a stationary measure}

An important consequence of the previous result is
the existence of a $P$-invariant Radon measure on the real line.
%
%th5.1 #&#
%
\begin{theorem}[(Existence of a $P$-invariant measure)]
\label{Cexistenceofinvariantmeasure}
Let $\mu$ be an irreducible, finitely supported, symmetric probability
measure on
$\mathrm{Homeo}^+({\mathbf R})$, and let $P$ be the associated transition
operator. Then there
exists a $P$-invariant Radon measure on the real line.
\end{theorem}
\begin{pf} Fix a continuous, compactly-supported function $\xi\dvtx
\mathbf R
\to[0,1]$ such that
$\xi\equiv1$ on $K$. For any initial point $x$, let us stop the
process $X^x_n$ at a
\textit{random} stopping time $T=T(w)$ chosen in a Markovian way so that,
for all $n \in\mathbb{N}$,
\[
\bbP(T=n+1 \mid T\ge n) = \xi\bigl(X^x_{n+1}\bigr).
\]
In other words, after each iteration of the initial random walk, when
we arrive at some point $y=X_{n+1}^x$
we stop with the probability $\xi(y)$, and we continue the iterations
with probability $1 - \xi(y)$.

Denote by $Y^x$ the random stopping point $X^x_T$, and consider its
distribution $p^x$
(notice that $T$ is almost surely finite since the process $X^x_n$
almost surely visits $K$ and
$\xi\equiv1$ on $K$). Due to the continuity of $\xi$, the measure
$p^x$ on $\mathbf R$
depends continuously (in the weak topology) on $x$. Therefore, the
corresponding diffusion
operator $P_{\xi}$ defined by
\[
P_{\xi}(\varphi) (x) = \E\bigl(\varphi\bigl(Y^x\bigr)\bigr)
= \int_{\mathbf R} \varphi(y) \,dp^x (y)
\]
acts on the space of continuous bounded functions on $\mathbf R$, and hence
it acts by duality on the
space of probability measures on $\mathbf R$. Notice that for any such
probability measure, its image
under $P_{\xi}$ is supported on $\widetilde{K}:= \supp(\xi)$. Thus,
applying the
Krylov--Bogolubov procedure of time averaging (and extracting a
convergent subsequence),
we see that there exists a $P_{\xi}$-invariant probability
measure $\nu_0$.

To construct a Radon measure that is stationary for the initial
process, we proceed as follows.
For each point $x\in\bbR$, let us take the sum of the Dirac measures
supported in its random
trajectory before the stopping\vspace*{1pt} moment $T$. In other words, we consider
the ``random measure''
$m_x (\omega):= \sum_{j=0}^{T(w)-1} \delta_{X_j^x}$. We then consider
its expectation
\[
m_x=\E\Biggl(\sum_{j=0}^{T(w)-1}
\delta_{X_j^x} \Biggr)
\]
as a measure on $\mathbf R$. Finally, we integrate $m_x$ with respect to
the measure $\nu_0$ on $x$,
thus yielding a Radon measure $\nu:= \int m_x \,d\nu_0(x)$ on $\mathbf
R$. Formally\vadjust{\goodbreak} speaking,
for any compactly supported function $f$, we have
%
%e1 #&#
%
\begin{equation}
\label{eqnudef} \int_{\mathbf R} f d \nu= \int_{\mathbf R}
\E\Biggl(\sum_{j=0}^{T(w)-1} f
\bigl(X_j^x\bigr) \Biggr) \,d\nu_0(x).
\end{equation}
Notice that the right-hand side expression of (\ref{eqnudef}) is well
defined and finite.
Indeed, there exist $N \in\mathbb{N}$ and $p_0 > 0$ such that with
probability at least
$p_0$ a trajectory starting at any point of $\supp(f)$ hits $K$ in at
most $N$ steps. Thus,
the distribution of the measure $m_x(w)$ on $\supp(f)$ [the number of
steps that are spent in
$\supp(f)$ until the stopping moment] has an exponentially decreasing
tail. Thus, its expectation
is finite and bounded uniformly on $x \in\supp(f)$, which implies the
finiteness of the integral.

Now, let us check that the measure $\nu$ is $P$-invariant. Let us
first rewrite
the measure $\nu$. To do this, notice that using the full probability formula,
one can check that the measure $m_x$ equals
\[
\sum_{n\ge0} \sum_{g_1,\ldots,g_n\in G}
\prod_{j=1}^n \mu(g_j) \cdot
\prod_{j=1}^{n} \bigl[ 1-\xi
\bigl(g_j\circ\cdots\circ g_1(x)\bigr) \bigr] \cdot
\delta_{g_n\circ\cdots\circ g_1(x)}.
\]
Thus
\begin{eqnarray*}
P (m_x) &=& \sum_{g \in G} \mu(g) \cdot
g_{\star} m_x
\\
&=& \sum_{g \in G} \mu(g) \cdot g_{\star} \Biggl(
\sum_{n\ge0} \sum_{g_1,\ldots,g_n\in G}
\prod_{j=1}^n \mu(g_j) \cdot
\prod_{j=1}^{n} \bigl[ 1-\xi
\bigl(g_j\circ\cdots\circ g_1(x)\bigr) \bigr] \\
&&\hspace*{237.2pt}{}\times
\delta_{g_n\circ\cdots\circ g_1(x)} \Biggr)
\\
&=& \sum_{n\ge0} \sum_{g_1,\ldots,g_n, g\in G}
\Biggl( \mu(g) \cdot\prod_{j=1}^n
\mu(g_j) \Biggr) \cdot\prod_{j=1}^{n}
\bigl[ 1-\xi\bigl(g_j\circ\cdots\circ g_1(x)\bigr) \bigr]\\
&&\hspace*{63.4pt}{}\times
g_{\star} \delta_{g_n\circ\cdots\circ g_1(x)}
\\
&=& \sum_{n\ge0}\sum_{g_1,\ldots,g_n,g_{n+1} \in G}
\Biggl( \prod_{j=1}^{n+1}
\mu(g_j) \Biggr) \cdot\prod_{j=1}^{(n+1)-1}
\bigl[ 1-\xi\bigl(g_j\circ\cdots\circ g_1(x)\bigr) \bigr]\\
&&\hspace*{75.4pt}{}\times
\delta_{g_{n+1}\circ g_n\circ\cdots\circ g_1(x)}.
\end{eqnarray*}
In the same way as before, one can check that the last expression
equals the expectation of
the random measure $\sum_{j=1}^{T(\omega)} \delta_{X_j^x}$. In this
sum, we are counting
the stopping time, but not the initial one, and therefore
\[
P m_x = m_x -\delta_x + \E(
\delta_{Y^x}).
\]
Integration with respect to $\nu_0$ yields
\begin{eqnarray*}
P \nu&=& P \biggl( \int_{\mathbf R} m_x \,d
\nu_0(x) \biggr) = \int_{\mathbf R} P(m_x) \,d
\nu_0(x)
\\[-2pt]
&=&\int_{\mathbf R} m_x \,d\nu_0(x) - \int
_{\mathbf R} \delta_x \,d\nu_0 (x) + \int
_{\mathbf R} \E(\delta_{Y^x}) \,d\nu_0(x) = \nu-
\nu_0 + P_{\xi} \nu_0.
\end{eqnarray*}
Since $\nu_0$ is $P_{\xi}$-invariant, we finally obtain $P\nu=\nu$,
as we wanted to show.~%
\end{pf}

%%%%%%%%%%%%%%%%%%%%%%%%%%%%%%%%%%%%%%%%%%%%%%%%%%%%%%%%%%%%%%%%%%%%%%%%%%%%%%%%%%%%%%%%%%%%%%%%%%%%%%%%%%%

%s6 #&#
\section{Properties of $P$-invariant measures}
\label{Spropertiesofinvariantmeasures}

This section is devoted to the study of the properties
of the $P$-invariant Radon
measures constructed in Section~\ref{Srecurrence}. The following
topological fact is well known;
we recall its proof for completeness.
%
%pr6.1 #&#
%
\begin{proposition}\label{Pminimal}
Let $G$ be a finitely generated, irreducible group of homeomorphisms of
the real line. Then
either $G$ carries a discrete orbit or there is a unique minimal
nonempty closed $G$-invariant
set $\mathcal M$. In the latter case, the closure of every orbit
contains $\mathcal M$.
\end{proposition}
\begin{pf}
Let $K$ be a compact interval that intersects every orbit; the existence
of such an interval follows from the proof of Theorem~\ref{Crecurrence}.
Let $\mathcal E$ be the family of nonempty compact subsets $\mathcal K$
of $K$ such that
$\mathcal K = (G \mathcal K) \cap K$. If $\{\mathcal{K}_\lambda,
\lambda\in\Lambda\}$
is a chain (with respect to inclusion) in $\mathcal{E}$, then
$\mathcal{K}_{\Lambda}:= \bigcap_{\lambda} \mathcal{K}_{\lambda
}$ also belongs
to $\mathcal{E}$. By Zorn's lemma, $\mathcal E$ has a maximal element
$\mathcal K_0$. Notice that $\mathcal M:= G \mathcal K_0$ is a nonempty
minimal $G$-invariant closed subset of the real line. When $\mathcal M$
is not a discrete set, every point of $\mathcal M$ is an
accumulation point, and $\mathcal M$ is locally compact. Hence,
there are only two possibilities when $G$ has no discrete
orbits---either $\mathcal M={\mathbf R}$
or $\mathcal M$ is locally homeomorphic to a Cantor set. In the first
case, the proposition
is proved. In the second case, the orbit of every point of $\mathcal
{M}$ is dense
in $\mathcal{M}$. We will now prove that the closure of orbits of
points in
${\mathbf R} \setminus\mathcal M$ contains $\mathcal M$; this will establish
the uniqueness of the set $\mathcal M$.
Let $C$ be an arbitrary connected component of ${\mathbf R} \setminus
\mathcal M$.
Then $C$ is bounded and its right endpoint $r$ belongs to $\mathcal M$.
Therefore,
there is a sequence of elements $g_n \in G$ such that $g_n (r)$ tends
to $r$ as
$n$ tends to infinity and $g_n (r) \neq r$ for every $n$ (otherwise,
the set of
accumulation points of the orbit of $r$ would be a closed $G$-invariant set
strictly contained in $\mathcal M$). Now $g_n(C)$ tends uniformly to
$r$ as $n$ tends to infinity. Since the closure of the orbit of $r$
equals~$\mathcal{M}$, this shows that the closure of the orbit
of any point in $C$ contains $\mathcal M$.
\end{pf}

If there is a discrete orbit, then the counting measure on it is a
Radon measure that is $G$-invariant, and in particular, it is also
$P$-invariant for any
probability measure $\mu$ on $G$. The next two lemmas provide converses
to this fact.\looseness=-1

%le6.2 #&#
%
\begin{lemma} \label{Latomicpart1}
Let $\mu$ be a symmetric probability measure
on $\mathrm{Homeo}^+ ({\mathbf R})$ whose support is finite and
generates an irreducible\vadjust{\goodbreak} group $G$. Let $\nu$ be a $P$-invariant Radon
measure on the real line. If there is a discrete orbit, then $\nu$
is supported on the union of discrete orbits and is totally
invariant.\vspace*{-2pt}
\end{lemma}
\begin{pf}
If there is a discrete orbit $O$, it can be parametrized by ${\mathbf
Z}$ and then the action of $G$ on $O$ is by integer translations. In
this situation, the normal subgroup $G^1$ formed by\vspace*{1pt} the
elements acting trivially on $O$ is recurrent by Polya's
theorem~\cite{Polya}. Let $\mu^1$ be the (symmetric) measure on $G^1$
obtained by balayage of $\mu$ to $G^1$. Observe that the restriction of
$\nu$ to each component $C$ of ${\mathbf R} \setminus O $ is a finite
measure that is invariant for the Markov chain induced by $\mu^1$ on
$C$. It now follows from Lemma \ref {Lbi-infiniteness} (or from
\cite{DKN}, Proposition 5.7) that this measure is supported on
$\operatorname{Fix}(G^1) \cap\ovverline{C}$, the set of global fixed
points for the group $G^1$ contained in the closure of $C$. As a
consequence, $G$ acts by ``integer translations'' on the support
of~$\nu$, which consists of discrete orbits. To see that $\nu$ is
invariant, notice that for each atom $x \in{\mathbf R}$, the function
$g \mapsto\nu(g(x))$ viewed as a function defined on $G / G^1
\sim{\mathbf Z}$ is harmonic and positive and hence constant.\vspace*{-2pt}
\end{pf}
%
%le6.3 #&#
%
\begin{lemma} \label{Latomicpart2}
Let $\mu$ be a symmetric probability measure on $\mathrm{Homeo}^+
({\mathbf R})$ whose support is finite
and generates an irreducible group $G$. Let $\nu$ be a $P$-invariant
Radon measure on the real line.
If the atomic part $\nu_a$ of $\nu$ is nontrivial, then it is supported
on a union of discrete orbits.\vspace*{-2pt}
\end{lemma}
\begin{pf}
Let $x\in{\mathbf R}$ be a point such that $\nu(x)>0$. Let $O = G(x)$ be
the orbit
of $x$ endowed with the discrete topology, and let $\ovverline{\nu}$
be the
measure on $O$ defined by $\ovverline{\nu}(y):= \nu(y)$. Then
$\ovverline
{\nu}$
is an invariant measure for the Markov process induced by $\mu$ on $O$.
%To see this, consider an arbitrary function $f\dvtx O\rightarrow{
%with finite support. Let $f_{\varepsilon}\dvtx{\mathbf R}\rightarrow{
%be a family
%of functions with support a union of intervals of length $\varepsilon$
%centered
%at the points of the support of $f$. Now,
%f_{\varepsilon} \,d\nu=
%which is sufficient to establish our claim.

Let $L$ be an arbitrary compact interval containing
the compact interval $K$ constructed in the proof of Theorem~\ref{Crecurrence}
and let $R:= L \cap O$. We want to show that $R$ is finite. To do
this, first observe
that $R$ is a recurrent subset of $O$, by Theorem~\ref{Crecurrence}.
Let $Y$ be the
Markov chain on $R$ defined by the first return of $X$ to $R$. This
Markov chain
is symmetric because $X$ is symmetric. Moreover, the restriction of
$\ovverline{\nu}$ to
$R$ is invariant. Now since $\sum_{y\in R} \ovverline{\nu}(y)<\infty$,
there must be an
atom $y\in R$ such that $\ovverline{\nu}(y) $ is maximal. The
$P_Y$-invariance of
$\ovverline{\nu}_{|R}$ and the symmetry of the transition probabilities
$p_Y(\cdot,\cdot)$ yield
\[
\sum_{z\in R} p_Y(y,z) \ovverline{\nu}(z) =
\sum_{z\in R} p_Y(z,y) \ovverline{\nu}(z) =
\ovverline{\nu}(y).
\]
The maximum principle now implies that $\ovverline{\nu}(z) = \ovverline
{\nu}(y)$ for all
$z \in O$. Thus, all the atoms of $\ovverline{\nu}$ contained in $R$
have the same mass and
hence there is only a finite number of them. In particular, this
argument shows that
$O$ is discrete.\vspace*{-2pt}
\end{pf}

Next we consider the case where $G$ has no discrete orbits. As in Proposition
\ref{Pminimal}, let $\mathcal{M}$ be the unique nonempty
minimal $G$-invariant closed subset of the real line.\vadjust{\goodbreak}
%
%le6.4 #&#
%
\begin{lemma} \label{Lsupport}
Let $\mu$ be an irreducible, symmetric measure on $\mathrm{Homeo}^+
({\mathbf R})$
whose support is finite and generates a group $G$ without discrete
orbits on the real line. Then any $P$-invariant Radon measure is
supported on $\mathcal M$.
\end{lemma}
\begin{pf}
Let $\nu$ be a $P$-invariant Radon measure on the real line.
The measure $\nu$ is quasi-invariant by $G$, because
for all $h$ in the support of $\mu$ we have
\[
h_{\star} \nu\leq\frac{1} {\mu(h)} \sum_{g \in G}
\mu(g) g_{\star}\nu= \frac{1}{\mu(h)}\nu.
\]
So the support of $\nu$ is a closed
$G$-invariant subset of the real line, and hence, it contains $\mathcal M$.
Therefore, it suffices to verify that $\nu$ does not charge any
component of
$\mathcal M^c$. If $\mathcal M^c$ is nonempty, we may collapse each connected
component of $\mathcal M^c$ to a point, thus obtaining a topological
real line
carrying a $G$-action for which every orbit is dense. The $P$-invariant
measure $\nu$
can be pushed to a $P$-invariant Radon measure $\ovverline{\nu}$ for
this new action. If
a component of $\mathcal M^c$ has a positive charge, then $\ovverline
{\nu}$ has atoms. By
Lemma~\ref{Latomicpart2}, this implies that the $G$-action cannot
be minimal after
the collapsing, which is a contradiction. We thus conclude that the
original $P$-invariant
measure $\nu$ does not charge the components of $\mathcal M^c$, and
so, $\nu$ must be
supported on $\mathcal M$.
\end{pf}

%%%%%%%%%%%%%%%%%%%%%%%%%%%%%%%%%%%%%%%%%%%%%%%%%%%%%%%%%%%%%%%%%%%%%%%%%%%%%%%%%%%%%%%%%%%%%%%%%%%%%%%

%s6.1 #&#
\subsection{Uniqueness of the $P$-invariant Radon measure} \label{Suniqueness}

When the action of $G$ possesses discrete orbits, we
know that
every stationary Radon measure must be $G$-invariant; however, two such measures
may be supported on different orbits. We now establish the uniqueness
(up to a
scalar factor) of the stationary measure in the case where $G$ is a finitely
generated, irreducible subgroup of $\mathrm{Homeo}^+({\mathbf R})$ without
discrete orbits. Recall that, in this case, there exists
a unique minimal closed $G$-invariant set $\mathcal M$, and the orbit
of every
point in $\mathcal M$ is dense in $\mathcal M$; see Proposition~\ref{Pminimal}.
%
%th6.5 #&#
%
\begin{theorem} \label{Tunicity}
Let $\mu$ be a symmetric measure on $\mathrm{Homeo}^+ ({\mathbf R})$
whose support
is finite and generates an irreducible group $G$ without discrete orbits.
Then the $P$-invariant Radon measure $\nu$ is unique up to a scalar factor,
and its support is $\mathcal M$. Moreover, for all continuous functions
$\varphi,\psi$ with compact support, with $\varphi\geq0$ and
$\int\varphi\,d\nu> 0$, and for every $x\in{\mathbf R}$,
we have a.s. the convergence
%
%e2 #&#
%
\begin{equation}
\label{Econvergence} \frac{\psi(X_1^x) + \cdots+ \psi(X_N^x)}{\varphi
(X_1^x)+ \cdots+
\varphi(X_N^x)} \longrightarrow\frac{\int\psi\,d\nu} {\int\varphi\,
d\nu}
\end{equation}
as $N$ tends to infinity.
\end{theorem}

For the proof of this theorem, we first consider the case when every
$G$-orbit is dense.
Let $\nu$ be a $P$-invariant measure. We know that $\nu$ is fully
supported and has no
atoms. By Lemma~\ref{Lbi-infiniteness},\vadjust{\goodbreak} we may consider the distance
$d$ on the
real line defined by
\[
d (x,y):= \nu[x,y],\qquad x \leq y.\vspace*{-2pt}
\]

%le6.6 #&#
%
\begin{lemma} \label{Lweakcontraction}
For any fixed number $0<p<1$ and all $x,y$,
with probability at least $p$ we have
\[
\lim_{n\rightarrow\infty} d \bigl( X_n^x, X_n^y
\bigr) \leq\frac{d (x,y)} {1-p}.\vspace*{-2pt}
\]
\end{lemma}
\begin{pf} Since $\nu$ is $P$-invariant, the sequence of random variables
$\omega\mapsto d(X_n^x,X_n^y)$ is a martingale. In particular, for
every integer
$n \geq1$ we have
\[
\mathbb E \bigl( d \bigl( X_n^x, X_n^y
\bigr) \bigr) = d(x,y).
\]
By the martingale convergence theorem, the sequence $d(X_n^x,X_n^y)$
converges a.s. to a nonnegative random variable $v(x,y)$.
By Fatou's inequality, for every $x<y$ we have
\[
\mathbb E \bigl(v (x,y)\bigr) \leq\lim_{n\rightarrow\infty} \mathbb E
\bigl( d
\bigl(X_n^x, X_n^y\bigr) \bigr) =
d (x,y).
\]
The lemma then follows from Chebyshev's inequality.\vspace*{-2pt}
\end{pf}

We will combine the preceding lemma with an argument from~\cite{Brofferio}.
For this, recall that a $P$-stationary measure $\nu$ is said to be
\textit{ergodic}
if every $G$-invariant measurable subset either has measure $0$ or its
complement
has measure $0$.\vspace*{-2pt}

%le6.7 #&#
%
\begin{lemma} Assume that the hypotheses of Theorem~\ref{Tunicity} are
satisfied and
that every $G$-orbit is dense. If $\nu$ is an ergodic $P$-invariant
Radon measure,
then the convergence (\ref{Econvergence}) holds a.s. for every $x
\in {\mathbf R}$.\vspace*{-2pt}
\end{lemma}
\begin{pf} The diffusion operator acting on $L^1(\mathbf R, \nu)$ is a
positive contraction. Moreover, because of the recurrence of the Markov process,
this operator is conservative. We may hence apply the Chacon--Ornstein theorem
\cite{Ch-Or}, which together with the ergodicity of $\nu$ shows that
for $\nu$-almost
every point $x\in{\mathbf R }$ and \textit{all} functions
$\varphi, \psi$ in $C_c ({\mathbf R})$ such that $\varphi\geq0$ and
$\varphi= 1 $ on the interval of recurrence $K$ constructed in
the proof of Theorem~\ref{Crecurrence}, we have almost surely
%
%e3 #&#
%
\begin{equation}
\label{eqequidistribution} \lim_{n\rightarrow\infty} \frac{S_n\psi
(x,\omega)} {
S_n \varphi(x,\omega)}=
\frac{\int\psi\,d\nu}{ \int\varphi
\,d\nu},
\end{equation}
where $S_n \psi(x,\omega):= \psi( X_1^x ) + \cdots+
\psi(X_n^x)$
[and similarly for $S_n \varphi(x,\omega)$]. Let $y \in
{\mathbf R}$ and
the functions $\varphi$, $\psi$ be fixed.\vadjust{\goodbreak} We claim that, for any $k
\geq1$,
with probability at least $1 - 1/k$ we have
%
%e4 #&#
%
\begin{equation}
\label{eqalmostequidistribution} \limsup_{n\rightarrow\infty} \biggl
\llvert
\frac{S_n \psi(y,\omega
)}{S_n\varphi(y,\omega) } - \frac{\int\psi\,d\nu}{\int\varphi\,d\nu}
\biggr\rrvert\leq\frac{1}{k}.
\end{equation}
This obviously implies that (\ref{eqequidistribution}) holds almost surely.\vadjust{\goodbreak}
%Before showing this estimate notice that 0--1-type arguments
%imply that (\ref{eqequidistribution}) holds with probability $1$
%for every $\delta$, thus establishing (\ref{Econvergence}) for $y
%from which the lemma easily follows.

Since $\nu$ has total support, one can find a point $x$ generic in
the sense of (\ref{eqalmostequidistribution}) and sufficiently
close to $y$ so that $d(x,y) \leq\varepsilon$. From
Lemma~\ref{Lweakcontraction}, with probability at least $1/2$ we have
for all $n$ sufficiently large, say $n \ge n_0(\omega)$,
%
%e5 #&#
%
\begin{equation}
\label{eqdist-3e} d\bigl(X_n^y,X_n^x
\bigr) \leq(k+1)\varepsilon.
\end{equation}
Now, as we already know that (with probability 1)
\[
\lim_{n\rightarrow\infty} \frac{S_n\psi(x,\omega)} {
S_n \varphi(x,\omega)}= \frac{\int\psi \,d\nu}{ \int\varphi
\,d\nu},
\]
instead of estimating the difference in (\ref{eqalmostequidistribution}),
it suffices to obtain estimates of the ``relative errors''
%
%e6 #&#
%
\begin{equation}
\label{eqrel1} \limsup_{n\rightarrow\infty} \biggl\llvert\frac{S_n \psi
(y,\omega) - S_n \psi(x,\omega)}{S_n\varphi
(x,\omega) } \biggr
\rrvert\leq\delta_1(\varepsilon)
\end{equation}
and
%
%e7 #&#
%
\begin{equation}
\label{eqrel2} \limsup_{n\rightarrow\infty} \biggl\llvert\frac{S_n
\varphi(y,\omega) - S_n \varphi(x,\omega
)}{S_n\varphi(x,\omega) } \biggr
\rrvert\leq\delta_2(\varepsilon)
\end{equation}
in such a way that $\delta_1 (\varepsilon) \to0$ and $\delta
_2(\varepsilon) \to0$
as $\varepsilon\to0$.

Since the estimate (\ref{eqrel2}) for $\varphi$ is a particular case
of the
estimate (\ref{eqrel1}), we will only check (\ref{eqrel1}).
Now, (\ref{eqdist-3e})
implies that
\begin{eqnarray*}
&&
\bigl\llvert S_n \psi(y,\omega) - S_n \psi(x,\omega)
\bigr\rrvert
\\
&&\qquad\le\operatorname{mod}\bigl((k+1)\varepsilon, \psi\bigr) \card
\bigl\{n_0(
\omega)\le j\le n \mid\mbox{either } X_j^x \mbox{ or }
X_j^y \mbox{ is in } \supp\psi\bigr\}\\
&&\qquad\quad{} + {2
n_0(\omega) \max}|\psi|
\\
&&\qquad\le\operatorname{mod}\bigl((k+1)\varepsilon, \psi\bigr) \card
\bigl\{j\le n \mid
X_j^x \in U_{(k+1)\varepsilon}(\supp\psi)\bigr\} + \const(
\omega).
\end{eqnarray*}
Here, $\operatorname{mod}(\cdot, \psi)$ stands for the modulus of continuity
of $\psi$ with respect
to the distance $d$ on the variable, and $U_{(k+1)\varepsilon}(\supp
\psi)$ denotes the
$(k+1)\varepsilon$-neighborhood of the support of $\psi$, again with
respect to $d$.

Let $\chi$ be a continuous function satisfying $0 \le\chi\le1$ and
that is equal to $1$ on
$U_{(k+1)\varepsilon}(\supp\psi)$ and to $0$ outside $U_{(k+2)
\varepsilon}(\supp\psi)$.
We have
\[
\card\bigl\{j\le n \mid X_j^x \in U_{(k+1)\varepsilon}(
\supp\psi)\bigr\} \le S_n \chi(x,\omega).
\]
Thus
\begin{eqnarray*}
&&
\biggl\llvert\frac{S_n \psi(y,\omega) - S_n \psi(x,\omega)}{S_n\varphi
(x,\omega) } \biggr\rrvert
\\
&&\qquad\le\frac{\const(\omega) + \operatorname{mod}((k+1)\varepsilon,
\psi)\cdot
S_n \chi(x,\omega)}{S_n\varphi(x,\omega)} \\
&&\qquad\xrightarrow{n\to\infty} {}
\operatorname{mod}\bigl((k+1)
\varepsilon, \psi\bigr)\cdot\frac{\int\chi \,d\nu
}{\int
\varphi \,d\nu} =: \delta_1(
\varepsilon).
\end{eqnarray*}
[Notice here that we have applied the fact that, by our choice of $x$, the
equality (\ref{eqequidistribution}) holds for the\vadjust{\goodbreak} functions $\chi$ and
$\varphi$.] Since $\operatorname{mod}((k+1)\varepsilon, \psi)$ tends to
$0$ as
$\varepsilon\to0$ and the quotient
\[
\frac{\int\chi \,d\nu}{\int\varphi \,d\nu} \le\frac{\nu
(U_{(k+2)\varepsilon}(\supp\varphi))}{\int\varphi
\,d\nu}
\]
remains bounded, this yields $\delta_1(\varepsilon)\to0$ as
$\varepsilon\to0$.
\end{pf}

It is now easy to finish the proof of Theorem~\ref{Tunicity} in the
case where all
the $G$-orbits are dense. Indeed, given any two ergodic $P$-invariant
Radon measures
$\nu_1,\nu_2$, for all $x \in{\mathbf R}$ and all compactly supported,
real-valued function $\psi$, we have almost surely
\[
\frac{S_N \psi(x,\omega)}{S_N \varphi(x,\omega)} \longrightarrow\frac
{\int\psi \,d\nu_i} {\int\varphi \,d\nu_i},
\]
where $i=1,2$. Thus, $\int\psi \,d\nu_1 = \lambda\int\psi \,d
\nu_2$,
with $\lambda:= \int\varphi \,d \nu_1 / \int\varphi \,d\nu_2$.
This proves that $\nu_1 = \lambda\nu_2$. The case
of nonnecessarily ergodic $\nu_1,\nu_2$
follows from standard ergodic decomposition type arguments.

The proof of Theorem~\ref{Tunicity} in the nonminimal case is more technical,
because the argument of collapsing the connected components of the
complement of the
unique minimal invariant closed set $\mathcal{M}$ is delicate. Indeed,
although this
procedure induces a minimal action from which the uniqueness of the
stationary measure
(up to a scalar factor) may be easily deduced, establishing (\ref
{Econvergence}) is
much more complicated, mainly due to the fact that, after collapsing,
the functions
$\psi,\varphi$ are no longer continuous. Below we propose two different
solutions
to this problem.
\begin{pf*}{First proof of Theorem~\ref{Tunicity} in the nonminimal case}
As before,
the main point consists in obtaining a good estimate of the form (\ref
{eqrel1}). To do this,
we fix $\varepsilon_0 > 0$, and we consider all the connected
components of the complement
of $\mathcal{M}$ over which the oscillation of $\psi$ is at least
$\varepsilon_0$. Since
$\psi$ has compact support, there are only finitely many such
components, say
$C_1,\ldots,C_k$. Given $\varepsilon_1 > 0$, let us consider a
continuous function
$\chi_1$ satisfying $0 \leq\chi_1 \leq1$ and that is equal to~$1$
on each
$U_{\varepsilon_1}(C_i)$ and to $0$ outside $\bigcup_i U_{2
\varepsilon
_1}(C_i)$.
Now, take $\varepsilon_2 > 0$ such that, if
$d (x,y) \leq3 \varepsilon_2$, then either
$|\psi(x)-\psi(y)| \leq2\varepsilon_0 \mbox{ or } x \mbox{
belongs to }
\bigcup_i U_{\varepsilon_1} (C_i)$.
(The existence of such an $\varepsilon_2$ is easy to establish.)
Finally, let $\chi$
be a continuous function satisfying $0 \leq\chi\leq1$ and that is
equal to $1$
on the set
\[
S_1:= \bigl\{x \mid\mbox{there is } y \in\supp\psi\mbox{ such that
} d (x,y) \leq3 \varepsilon_2 \bigr\}
\]
and to $0$ outside
$\{x \mid d (x,y) \geq4 \varepsilon_2 \mbox{ for all } y \in\supp
\psi\}$.

Notice that, although $d$ is not a metric on the line, we still have
that, if
\mbox{$d (x,y) \leq\varepsilon$}, then with probability at least $1-1/k$
there is
$n_0(\omega)$ such that, for all $n \geq n_0(\omega)$,
\[
d \bigl(X_n^x,X_n^y\bigr)
\leq(k+1) \varepsilon.\vadjust{\goodbreak}
\]
Fix $\varepsilon\leq3\varepsilon_2/(k+1)$. Given
$y \in\mathbf R$, take a point $x$ that is generic in the sense of
(\ref{eqequidistribution}) and such that $d(x,y) \leq\varepsilon$.
With probability at least $1 - 1/k$ we have
\begin{eqnarray*}
\bigl|S_n \psi(x,\omega) - S_n \psi(y,\omega)\bigr| &\leq&\sum
_{j=1}^{n} \bigl|\psi\bigl(X_j^x
\bigr) - \psi\bigl(X_j^y\bigr)\bigr|
\\
&\leq&\operatorname{const}(\omega) + 2 \max|\psi| \operatorname{card}
\biggl\{ j\leq n \Bigm|
X_j^x \in\bigcup_i
U_{\varepsilon_1}(C_i) \biggr\}\\
&&{} + 2 \varepsilon_0
\operatorname{card} \bigl\{ j\leq n \mid X_j^x \in
S_1 \bigr\}
\\
&\leq&\operatorname{const}(\omega) + 2 \max|\psi| S_n \chi_1 (x,
\omega) + 2\varepsilon_0 S_n \chi(x,\omega).
\end{eqnarray*}
Dividing by $S_n \varphi(x,\omega)$ and passing to the limit we obtain
\begin{eqnarray*}
\limsup_{n \to\infty} \frac{|S_n \psi(x,\omega) - S_n \psi(y,\omega
)|}{S_n \varphi
(x,\omega)} &\leq& 2 \max|\psi| \frac{\int\chi_1 \,d\nu}{\int
\varphi\,
d \nu}
+ 2 \varepsilon_0 \frac{\int\chi \,d\nu}{\int\varphi
\,d\nu
}
\\
&\leq& \frac{2 \max|\psi|}{\int\varphi \,d\nu} \sum_i \nu\bigl(
U_{2\varepsilon_1}(C_i) \bigr) + 2 \varepsilon_0
\frac{\int\chi \,d\nu}{\int\varphi
\,d\nu}.
\end{eqnarray*}
To conclude, notice that the first term can be made arbitrarily small by
taking $\varepsilon_1$ very small, since the $\nu$-measure of the set
$\bigcup_i C_i$ is zero.
\end{pf*}
\begin{pf*}{Second proof of Theorem~\ref{Tunicity} in the
nonminimal case}
To obtain an estimate of the form (\ref{eqrel1}), we collapse the connected
components of $\mathcal{M}^c$, thus obtaining a topological real line carring
a minimal $G$-action. However, after collapsing, the functions $\psi$ and
$\varphi$ are no longer continuous. To solve this problem, we consider a
nonnegative function $\varphi_1 \in C_c (\mathbf R)$ that is positive
on the
recurrence interval $K$ and is contant on each connected component of
$\mathcal{M}^c$. If we are able to estimate (\ref{eqrel1}) but for
$\varphi_1$ instead of $\varphi$ and for any function $\psi$,
then we will have
\[
\frac{S_n \psi(y,\omega)}{S_n \varphi(y,\omega)} = \frac{S_n \psi
(y,\omega) / S_n \varphi_1 (y,\omega)} {
S_n \varphi(y,\omega) / S_n \varphi_1 (y,\omega)} \longrightarrow\frac
{\int\psi \,d\nu/ \int\varphi_1 \,d\nu} {
\int\varphi \,d\nu/ \int\varphi_1 \,d\nu} =
\frac{\int\psi \,d\nu}{\int\varphi \,d\nu}
\]
as we want to show.

Fix $\varepsilon> 0$. We leave to the reader the task of showing the existence
of $\psi_1,\chi_1$ in $C_c (\mathbf R)$ that are constant on each connected
component
of $\mathcal{M}^c$ and satisfy:
\begin{longlist}[(ii)]
\item[(i)]
$|\psi- \psi_1| \leq\chi_1$;

\item[(ii)] $\int\chi_1 \,d\nu\leq\varepsilon$.
\end{longlist}
Then using
\[
\frac{S_n \psi(y,\omega)}{S_n \varphi_1 (y,\omega)} = \frac{S_n \psi_1
(y,\omega)}{S_n \varphi_1 (y,\omega)} + \frac{S_n (\psi- \psi
_1)(y,\omega)}{S_n \varphi_1 (y,\omega)}\vadjust{\goodbreak}
\]
we obtain
\begin{eqnarray*}
&&
\biggl\llvert\frac{S_n \psi(y,\omega)}{S_n \varphi_1 (y,\omega)} -
\frac{\int\psi \,d\nu}{\int\varphi_1 \,d\nu} \biggr\rrvert\\
&&\qquad\leq
\biggl\llvert\frac{S_n\psi_1(y,\omega)}{S_n \varphi_1 (y,\omega)} -
\frac{\int\psi_1 \,d\nu}{\int\varphi_1 \,d\nu} \biggr\rrvert+ \biggl
\llvert\frac{S_n(\psi-\psi_1)(y,\omega)}{S_n\varphi_1(y,\omega)} \biggr
\rrvert+ \biggl\llvert\frac{\int(\psi_1 - \psi) \,d\nu}{\int\varphi_1
\,d\nu} \biggr
\rrvert
\\
&&\qquad\leq \biggl\llvert\frac{S_n\psi_1(y,\omega)}{S_n \varphi_1 (y,\omega
)} - \frac{\int\psi_1 \,d\nu}{\int\varphi_1 \,d\nu} \biggr\rrvert+
\biggl\llvert\frac{S_n \chi_1 (y,\omega)}{S_n \varphi_1 (y,\omega)} -
\frac{\int\chi_1 \,d\nu}{\int\varphi_1 \,d\nu} \biggr\rrvert+ 2 \biggl
\llvert\frac{\int\chi_1 \,d\nu}{\int\varphi_1 \,d\nu} \biggr\rrvert.
\end{eqnarray*}
As in the minimal case, we have
\[
\biggl\llvert\frac{S_n\psi_1(y,\omega)}{S_n \varphi_1 (y,\omega)} -
\frac{\int\psi_1 \,d\nu}{\int\varphi_1 \,d\nu} \biggr\rrvert
\longrightarrow0,\qquad \biggl\llvert\frac{S_n \chi_1 (y,\omega)}{S_n \varphi
_1 (y,\omega)} - \frac{\int\chi_1 \,d\nu}{\int\varphi_1 \,d\nu} \biggr
\rrvert\longrightarrow0.
\]
Since
\[
\biggl\llvert\frac{\int\chi_1 \,d\nu}{\int\varphi_1 \,d\nu} \biggr
\rrvert\leq\frac{\varepsilon}{\int\varphi_1 \,d\nu},
\]
this shows the desired convergence
\[
\frac{S_n \psi(y,\omega)}{S_n \varphi_1 (y,\omega)} \longrightarrow\frac
{\int\psi \,d\nu}{\int\varphi_1 \,d\nu},
\]
thus finishing the proof.
\end{pf*}

%%%%%%%%%%%%%%%%%%%%%%%%%%%%%%%%%%%%%%%%%%%%%%%%%%%%%%%%%%%%%%%%%%%%%%%%%%%%%%%%%%%%%%%%%%%5

%s7 #&#
\section{Global stability at a finite distance}
\label{Sstability}

We say that an irreducible subgroup $G$ of $\mathrm
{Homeo}^+({\mathbf R})$
has the \textit{strong contraction property} if there exists a compact
interval $L$
such that, for every compact interval $I$, there is a sequence of
elements $h_n$ of
$G$ such that $h_n(I) \subset L$ for all $n$, and the diameter of $h_n(I)$
tends to zero as $n$ tends to infinity. The group $G$ has the weak
\textit{contraction property} if the property above holds for all
compact intervals $I$ of length less than $1$.

For example, every non-Abelian subgroup of the affine group has the
strong contraction property. In the opposite direction,
the group $\widetilde{\mathrm{Homeo}}{}^+({\mathbf S}^1)$ of
homeomorphisms of the real line commuting with
the translation $x \mapsto x+1$ does not have the strong
contraction property, since no interval of length greater than $1$ can be
contracted to an interval of length less than $1$.
However, this group has the weak contraction property.

Recall that the action of a subgroup $G \subset\mathrm
{Homeo}^+({\mathbf R})$
is \textit{semi-conjugate} to that of an homomorphic image
$\ovverline{G} \subset\mathrm{Homeo}^+ ({\mathbf R})$ if there
exists a surjective, nondecreasing, continuous map $D\dvtx {\mathbf R}
\to
{\mathbf R}$
such that $D (g(x)) = \bar{g}(D(x))$ for all $x \in{\mathbf R}$ and all
$g \in G$, where $\bar{g}$ denotes the image of $g$ under the homomorphism.
(We have already met this situation in the proof of Lemma~\ref{Lsupport}.)
The following result was obtained by Malyutin~\cite{Malyutin}, although
an analogous
statement due to McCleary (see, e.g.,~\cite{glass}, Theorem 7.E) was
already known in the context of orderable groups. We include a proof for
completeness.
%
%th7.1 #&#
%
\begin{theorem}\label{Ttopologicalstructure}
Let $G$ be a finitely generated, irreducible
subgroup of $\mathrm{Homeo}^+ ({\mathbf R})$. Then
one of the following possibilities occur:
\begin{itemize}
\item$G$ has a discrete orbit;

\item$G$ is semi-conjugate to a minimal group of translations;

\item$G$ is semi-conjugate to a subgroup of
$\widetilde{\mathrm{Homeo}}{}^+({\mathbf S}^1)$ having the weak contraction
property;

\item$G$ has the strong contraction property.
\end{itemize}
\end{theorem}
\begin{pf} Assume that there are no discrete orbits.
By Proposition~\ref{Pminimal}, there is a unique minimal nonempty closed
$G$-invariant subset $\mathcal M$. Now collapse\vspace*{1pt} each connected component
of $\mathcal M^c$ to a point to semi-conjugate $G$ to a group
$\ovverline{G}$ whose action is minimal.
If $G$ preserves a Radon measure, then after semi-conjugacy this measure
becomes a $\ovverline{G}$-invariant Radon measure of total support and no
atoms. Therefore, $\ovverline{G}$ (resp., $G$) is conjugate (resp.,
semi-conjugate)
to a group of translations.

Now suppose that $G$ has no invariant Radon measure. We claim that the action
of $\ovverline{G}$ cannot be free. If the action was free, $\ovverline
{G}$ would be conjugate to a
group of translations by H\"older's theorem; see either~\cite{ghys} or
\cite{navas}.
Pulling back the Lebesgue measure by the semi-conjugacy would provide a
$G$-invariant
Radon measure, which is contrary to our assumption. So the action of
$\ovverline{G}$ is not free.

Let $\bar{g} \in\ovverline{G}$ be a nontrivial element having fixed points,
and let $\bar{x}_0$ be a point in the boundary of $\operatorname
{Fix}(\bar
{g})$. Then
there is a left or right neighborhood $I$ of $\bar{x}_0$ that is contracted
to $\bar{x}_0$ under iterates of either $\bar{g}$ or its inverse. By
minimality, every $\bar{x}$ has a neighborhood that can be contracted to
a point by elements in $\ovverline{G}$. Coming back to the original action,
we conclude that every $x \in{\mathbf R}$ has a neighborhood that can be
contracted to a point by elements in~$G$. Since $G$ is finitely generated,
such a point can be chosen to belong to a compact interval $L$ that
intersects every orbit (compare with Theorem~\ref{Crecurrence}).

For each $x\in{\mathbf R}$ define $T(x) \in{\mathbf R} \cup\{+\infty\}$
as the supremum of the $y>x$ such that the interval $(x,y)$ can be
contracted to a point in $L$ by elements of~$G$. Then either
$T \equiv+\infty$, in which case the group $G$ has
the strong contraction property, or $T(x)$ is finite for every $x\in
{\mathbf R}$.
In the last case, $T$ induces a nondecreasing map $\bar{T}\dvtx
{\mathbf R}
\rightarrow{\mathbf R}$
commuting with all the elements in $\ovverline{G}$. Since the union of
the intervals on which
$\bar{T}$ is constant is invariant by $\ovverline G$, the minimality of
the action implies
that there is no such interval, that is, $\bar T$~is strictly
increasing. Moreover,\vspace*{1pt}
the interior of ${\mathbf R} \setminus\bar{T} ({\mathbf R})$ is also
invariant, hence
empty because the action is minimal. In other words, $\bar T$ is continuous.
All of this shows that $\bar T$ induces a homeomorphism of ${\mathbf R}$
into its
image. Since the image of $\bar T$ is $\ovverline G$-invariant, it must
be the whole line.
Therefore, $\bar T$~is a homeomorphism from the real line to itself.
Observe that
$\bar{T} (x)>x$ for any point $x$, which implies that $\bar T$ is
conjugate to
the translation $x \mapsto x+1$. After this conjugacy, $\ovverline{G}$
becomes a
subgroup of $\widetilde{\mathrm{Homeo}}{}^+({\mathbf S}^1)$. This completes
the proof.
\end{pf}

We now establish a probabilistic version of Theorem~\ref{Ttopologicalstructure}.
Notice that in the first two cases given by this theorem, the Markov
chain $X$
induces a random walk on a (finitely generated) subgroup of the group
of translations.
In the other two cases, we establish the global stability at a finite distance.
More precisely, we obtain the following result.
%
%th7.2 #&#
%
\begin{theorem} \label{Tprobabilisticcontractionproperty}
Let $\mu$ be an irreducible, finitely supported, symmetric probability
measure on
$\mathrm{Homeo}^+({\mathbf R})$ such that the group $G$ generated by the
support of $\mu$ acts
minimally on ${\mathbf R}$. If $G$ has the strong contraction property,
then for any $x<y$ and
any compact interval $J$, almost surely we have
%
%e8 #&#
%
\begin{equation}
\label{Econtraction} {\mathbf1}_{J} \bigl(X_n^x
\bigr) \bigl| X_n^y - X_n^x \bigr|
\longrightarrow0 \qquad\mbox{as } n \to\infty.
\end{equation}
If $G$ satisfies only the weak contraction property, then viewed (after
conjugacy)
as a subgroup of $\widetilde{\mathrm{Homeo}}{}^+({\mathbf S}^1)$,
convergence (\ref{Econtraction}) holds with positive probability for
any $x < y < x+1$.
\end{theorem}

We will assume below that $G$ has the strong contraction property,
since the case
of the weak contraction property is analogous and may be left to the
reader. Moreover,
the result in the latter context is not new. Indeed, by conjugacy into
a subgroup of
${\widetilde{\mathrm{Homeo}}}{}^+({\mathbf S}^1)$, the $P$-invariant Radon
measures become invariant
by the translation $x \mapsto x+1$. Therefore, these mesures are
proportional to the ``pull-back''
of the unique $P$-stationary probability measure of the associated
action of $G$ on the
circle ${\mathbf R}/{\mathbf Z}$. Furthermore, for this associated
action, a
natural property of strong contraction for random compositions
holds. See~\cite{DKN}, Section 5.1, for more details.

The main technical ingredient of the proof of Theorem
\ref{Tprobabilisticcontractionproperty} is the next lemma,
which has an obvious extension to more general Markov processes.
%
%le7.3 #&#
%
\begin{lemma} \label{inf-often+good-exit} In the context of Theorem
\ref{Tprobabilisticcontractionproperty}, assume that $G$ has the
strong contraction
property, and let $K$ be any compact interval of recurrence. Fix $k \in
\mathbb{N}$ and
$h_1,\ldots,h_k$ in the support of $\mu$. Then almost surely the
following happens for
infinitely many $n \geq0$: the point $X_n^x$ belongs to $K$ and
$g_{n+1},\ldots,g_{n+k}$
coincide with $h_1,\ldots,h_k$, respectively.
\end{lemma}
\begin{pf}
Due to the Markov property, it suffices to show that this
situation almost surely happens at least once. Let $\xi\dvtx \mathbf{R}
\to[0,1]$ be the function defined by letting $\xi(z)$ be the
probability that there exists $n \geq0$ such that $X_n^z \in K$ and
$g_{n+i} = h_i$ for $i=1,\ldots,k$. We need to show that $\xi(x) = 1$,
and we will actually show that $\xi(z) = 1$ holds for all $z
\in\mathbf{R}$. To do this, let $p:= \mu(h_1) \cdots\mu(h_k) > 0$. For each
$\omega\in\Omega$ and $z \in\mathbf{R}$,\vadjust{\goodbreak} let $n(z) \geq0$ be the
first-entry time of $z$ into $K$. A moment reflexion shows that
\[
\xi(z) = p + (1 - p) \mathbb{E} \bigl( \xi\bigl(X_{n(z)+k}^z
\bigr) | (g_{n(z)+1},\ldots,g_{n(z)+k} ) \neq(h_1,\ldots,h_k) \bigr).
\]
Letting $\Phi:= \inf_{z \in\mathbf{R}} \xi(z)$, this yields
\[
\Phi\geq p + (1-p) \Phi.
\]
This easily implies that $\Phi=1$, as we wanted to show.
\end{pf}

The proof of Theorem~\ref{Tprobabilisticcontractionproperty} when
$G$ has the strong
contraction property is now easy. Indeed, let $L$ be the interval of
contraction,
and let $K:=[a,b]$ be a compact interval of recurrence containing $J$
and $x,y$.
As in Lemma~\ref{Lweakcontraction}, the value of $\nu(f_n (K))$ converges
to a limit $l(\omega)$ almost surely. Given $\omega\in\Omega$ for
which it
converges, we fix $M>0$ such that $\nu(f_n (K)) \leq M$ holds for all
$n$. Choose an interval $I:=[\bar{a},\bar{b}]$ containing $K$ so that
$\nu[\bar{a},a] > M$ and $\nu[b,\bar{b}] > M$. Given
$\varepsilon> 0$, let $h \in\Gamma$ be such that $h(I) \subset J$ and
$|h(I)| \leq\varepsilon$. Write $h$ in the form $h_k \cdots h_1$,
where each
$h_i$ is in the support of $\mu$. Lemma~\ref{inf-often+good-exit} shows
that if $\omega$ is generic there exist infinitely many $n \in
\mathbb{N}$ such that:
\begin{itemize}
\item$f_n (x)$ belongs to $K$;
\item$\nu(f_{n} (K)) \leq M$;
\item$(g_{n+1},\ldots,g_{n+k}) = (h_1,\ldots,h_k)$.
\end{itemize}
Since $f_n(K)$ intersects $K$ and its $\nu$-measure
is bounded from above by $M$, it must be contained in $I$. Therefore,
\[
f_{n+k} (K) = g_{n+k} \cdots g_{n+1} f_n
(K) = h_k \cdots h_1 f_n(K) \subset
h_{k} \cdots h_{1} (I) = h(I) \subset L,
\]
hence
\[
\bigl| f_{n+k} (K) \bigr| \leq\bigl| h (I) \bigr| < \varepsilon.
\]
Since $n$ can be taken as large as required and $\nu$ has no atoms, we must
necessarily have $l(\omega) \leq\varepsilon$. Since this is true for all
$\varepsilon> 0$, we conclude that $l(\omega) = 0$. This
implies the desired result.

It should be emphasized that the distance $d$ induced by $\nu$ and the usual
distance on ${\mathbf R}$ may be very different in neighborhoods of $\pm
\infty$. As an example,
consider the case of a non-Abelian subgroup $G$ of the affine group
generated by
an expansion and a translation. Writing $g(x) = a_g x + b_g$ for each
$g \in G$,
the homomorphism $g \mapsto\log(a_g)$ induces a (symmetric) random
walk on ${\mathbf Z}$,
which is therefore recurrent. As a consequence, the length of the interval
$[X_n^x,X_n^y]$ oscillates between $0$ and $\infty$ even though its
$\nu
$-measure
converges to zero.

%%%%%%%%%%%%%%%%%%%%%%%%%%%%%%%%%%%%%%%%%%%%%%%%%%%%%%%%%%%%%%%%%%%%%%%%%%%%%%%%%%%%%%%%%%%%%%%%%%%%%%

%s8 #&#
\section{Derriennic's property and Lipschitz actions}
\label{Sderriennic}

Let $\mu$ be a symmetric probability measure on
$\mathrm
{Homeo}^+({\mathbf R})$
with finite support generating an irreducible group~$G$. We will say
that the pair $(G,\mu)$\vadjust{\goodbreak}
has the \textit{Derriennic property} if, for every $x\in{\mathbf R}$,
\[
x= \int_{G} g(x) \,d\mu(g).
\]
This terminology is inspired by~\cite{Derriennic}, where Derriennic
studies Markov
processes on the real line satisfying $\mathbb E ( X^x_1 ) = x$ for
large values
of $|x|$. As we demonstrate below, under very general conditions, this property
is always guaranteed after a suitable semi-conjugacy.
%
%pr8.1 #&#
%
\begin{proposition}\label{Tderrienniccoordinates}
Let $\mu$ be a finitely supported, symmetric measure on $\mathrm
{Homeo}^+ ({\mathbf R})$
whose support generates an irreducible group $G$ without discrete
orbits. Then $G$ is
semi-conjugate to a group $\ovverline{G}$ so that the pair $(\ovverline
{G},\mu)$ has
the Derriennic property.
\end{proposition}
\begin{pf} Since $G$ has no discrete orbits, there is a unique nonempty
$G$-invariant
closed subset $\mathcal M$ in which every orbit is dense. By Lemma \ref
{Lsupport},
the support of the $P$-invariant measure $\nu$ coincides with
$\mathcal M$.
Moreover, Lemma~\ref{Latomicpart2} shows that no $P$-invariant measure
$\nu$ has atoms. Now fix a point $x_0$ in the real line and consider
the map
\[
x \in{\mathbf R} \mapsto D(x):= \cases{ \nu[x_0,x], &\quad if $x \geq
x_0$,
\cr
-\nu[x,x_0], &\quad if $x \leq x_0$. }
\]
This map is continuous
and nondecreasing. Furthermore, Lemma~\ref{Lbi-infiniteness} implies
that this map is
also surjective. Consequently, since the support of $\nu$ is
$G$-invariant, $D$ induces a
semi-conjugacy from $G$ to a group $\ovverline{G}$ whose action is
minimal. We claim that
the pair $(\ovverline{G},\mu)$ has the Derriennic's property. Let
$\ovverline{P}_{\mu}$ be
the transition operator associated to the Markov process. Notice that
$D$ maps the measure
$\nu$ to the Lebesgue measure, which is then $\ovverline{P}_{\mu
}$-invariant. Now, for any
$x<y$, we have
\[
y-x = \int_{\ovverline G} \bigl( \bar{g} (y) - \bar{g} (x) \bigr) \,d
\mu(\bar{g}),
\]
which implies that the value of the \textit{drift},
\[
\operatorname{Dr} (\ovverline{G},\mu):= \int_{\ovverline G} \bigl( {\bar g}(x) -x
\bigr) \,d\mu(\bar{g})
\]
is independent of $x$. To complete the proof, we need to show that the
drift vanishes. To do this, we closely follow the argument of the first
proof of~\cite{DKN}, Proposition 5.7.

Fix any $a<b$, and let us integrate (\ref{eqdrift}) over $[a,b]$, then
doubling the integral
in order to couple $g$ and $g^{-1}$,
\begin{eqnarray*}
2 \int_a^b \operatorname{Dr}(\bar{G}, \mu) \,dx &=& \int
_a^b \biggl( \int_G
\bigl(\bar{g}(x)-x\bigr) \,d\mu(\bar{g}) + \int_G \bigl(
\bar{g}^{-1}(x)-x\bigr) \,d\mu(\bar{g}) \biggr) \,dx
\\
&=& \int_G \biggl( \int_a^b
\bigl[ \bigl(\bar{g}(x)-x\bigr)+ \bigl(\bar{g}^{-1}(x)-x\bigr) \bigr]
\,dx \biggr) \,d\mu(\bar{g}).
\end{eqnarray*}
Now, we will transform the value under the integral by means of the
following notion.
%
%de8.2 #&#
%
\begin{definition}
For any $c\in\R$ and $\bg\in\mathrm{Homeo}^+({\mathbf R})$, let
\begin{eqnarray*}
\Phi_{\bg}(c)&=&\Phi_{\bg^{-1}}(c) \\
&=& \mes\bigl\{ (x,y) \mid\mbox{either }
x<c<y<\bg(x) \mbox{ or } x<c<y<\bg^{-1}(x) \bigr\},
\end{eqnarray*}
where $\mes$ is the two-dimensional Lebesgue measure; see Figure~\ref{figPhi}.
Equivalently,
\[
\Phi_{\bg}(c)=\Phi_{\bg^{-1}}(c) = \cases{ \displaystyle\int
_{\bar{g}^{-1}(x)}^{x} \bigl[ \bar{g}(s) - s \bigr] \,ds, &\quad if $
\bar{g} (x) \geq x$,
\vspace*{2pt}\cr
\displaystyle\int_{\bar{g}(x)}^x \bigl[
\bar{g}^{-1}(s) - s \bigr] \,ds, &\quad if $\bar{g}(x) \leq x$. }
\]
\end{definition}

A geometric argument based on symmetry yields the following lemma.
%
%le8.3 #&#
%
\begin{lemma}\label{lPhi-int}
For any $\bg\in\mathrm{Homeo}^+({\mathbf R})$ and for any interval $[a,b]$
we have
%
%e9 #&#
%
\begin{equation}
\label{eqdrift} \int_a^b \bigl[ \bigl(
\bar{g}(x)-x\bigr)+ \bigl(\bar{g}^{-1}(x)-x\bigr) \bigr] \,dx =
\Phi_{\bg}(b)- \Phi_{\bg}(a).
\end{equation}
\end{lemma}
\begin{pf}
Notice that $\int_a^b (\bg(x)-x) \,dx$ equals
\[
\mes\bigl\{(x,y) \mid a<x<b, x<y<\bg(x) \bigr\} - \mes\bigl\{ (x,y)
\mid
a<x<b, \bg(x)<y<x \bigr\},
\]
which may be rewritten as
%
%e10 #&#
%
\begin{eqnarray}
\label{eqintrep}
&&\mes\bigl\{(x,y) \mid a < x < b, b < y < \bg(x) \bigr\} \nonumber\\
&&\qquad{}+ \mes
\bigl\{ (x,y) \mid a < x < b, a < y < b, x < y <\bg(x) \bigr\}
\nonumber\\[-8pt]\\[-8pt]
&&\qquad{}- \mes\bigl\{ (x,y) \mid a < x < b, \bg(x) < y < a \bigr\} \nonumber\\
&&\qquad{}-
\mes\bigl\{ (x,y)
\mid a < x < b, a < y < b, \bg(x) < y < x \bigr\}.
\nonumber
\end{eqnarray}
A similar
equality holds when changing $\bg$ by $\bg^{-1}$. Now, when taking the
sum of
$\int_a^b (\bg(x)-x) \,dx$ and $\int_a^b (\bg^{-1}(x)-x) \,dx$, we
see that the corresponding
second and fourth terms from (\ref{eqintrep}) cancel each other.
Indeed, these terms correspond
to the couples $(x,y) \in[a,b]^2$, and we have $x<y<\bg(x)$ if and
only if $\bg^{-1}(y)<x<y$.
The symmetry argument then shows that the second term for $\bg$ is
exactly the negative of
the fourth term for $\bg^{-1}$, and vice versa.

%f1 #&#
%
\begin{figure}

\includegraphics{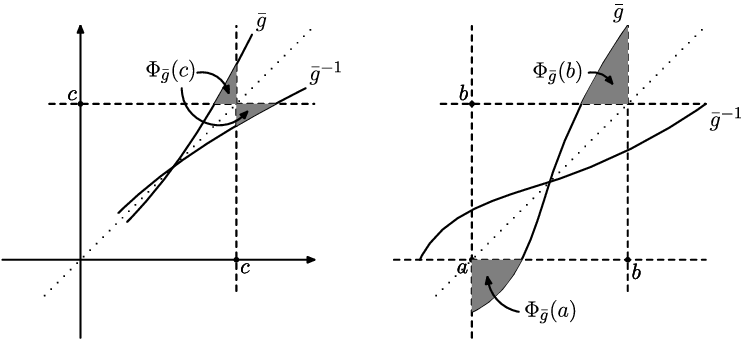}

\caption{Definition of $\Phi_{\bg}$ and illustration for the proof of
Lemma \protect\ref{lPhi-int}.}\label{figPhi}
\end{figure}

Therefore, the value of
\[
\int_a^b \bigl[ \bigl(\bar{g}(x)-x\bigr)+
\bigl(\bar{g}^{-1}(x)-x\bigr) \bigr] \,dx
\]
equals
\begin{eqnarray*}
&& \bigl[ \mes\bigl\{ (x,y) \mid a < x < b, b < y < \bg(x) \bigr\}\\
&&\qquad\hspace*{0pt}{} +
\mes\bigl
\{ (x,y) \mid a < x < b, b < y < \bg^{-1}(x) \bigr\} \bigr]
\\
&&\qquad{}- \bigl[ \mes\bigl\{ (x,y) \mid a < x < b, \bg(x) < y < a
\bigr\} \\
&&\qquad\hspace*{15.1pt}{}+ \mes
\bigl\{ (x,y) \mid a < x < b, \bg^{-1}(x) < y < a \bigr\} \bigr],
\end{eqnarray*}
and one can easily see that the expressions inside the brackets are equal
to $\Phi_{\bg}(b)$ and $\Phi_{\bg}(a)$, respectively; see
Figure~\ref{figPhi}.
This proves the desired equality.
\end{pf}

Now, we can complete the proof of Proposition~\ref{Tderrienniccoordinates}:
by integrating (\ref{eqdrift}) over $G$ we obtain, for any $a<b$,
\[
2(b-a) \operatorname{Dr}(\bar{G}, \mu) = \int_G \bigl(
\Phi_{\bg}(b) - \Phi_{\bg}(a)\bigr) \,d\mu(\bg).
\]
Denoting now $\Phi_{\mu}(c):= \int_G \Phi_{\bg}(c) \,d\mu(\bg
)$, this yields
\[
2(b-a) \operatorname{Dr}(\bar{G}, \mu) = \Phi_{\mu}(b) - \Phi_{\mu}(a).
\]
The last equality shows that $\Phi_{\mu}$ is an affine function. On the
other hand,
$\Phi_{\mu}$ is an average of nonnegative functions, and thus it is
nonnegative.
Therefore, $\Phi_{\mu}$ must be constant, which implies that $\operatorname{Dr}(\bar
{G},\mu)=0$.
\end{pf}

%%%%%%%%%%%%%%%%%%%%%%%%%%%%%%%%%%%%%%%%%%%%%%%%%%%%%%%%%%%%%%%%%%%%%%%%%%%%%%%%%%%%%%%%%%%%%%%%%

%%%%
%%%% THE PICTURES
%%%%

%%%%%%%%%%%%%%%%%%%%%%%%%%%%%%%%%%%%%%%%%%%%%%%%%%%%%%%%%%%%%%%%%%%%%%%%%%%%%%%%%

The next proposition demonstrates the relevance of the Derriennic property
in the study of the smoothness of a group action.
%
%pr8.4 #&#
%
\begin{proposition} \label{Llipschitz}
If a pair $(G,\mu)$ has the Derriennic property, then every element of $G$
is a Lipschitz map. Moreover, the displacement function $x \mapsto g(x)
- x$
is uniformly bounded in $x$ for every $g \in G$.\vadjust{\goodbreak}
\end{proposition}
\begin{pf} It suffices to prove the lemma for the elements of the
support of $\mu$.

To check the Lipschitz property, notice that
for any $g_0\in\supp\mu$ and any $x<y$ we have
\[
\mu(g_0) \cdot\bigl(g_0(y)-g_0(x)\bigr)
\le\int_G \mu(g) \cdot\bigl(g(y)-g(x)\bigr) \,d\mu(g) = y-x
\]
and thus
\[
g_0(y)-g_0(x) \le\frac{1}{\mu(g_0)} \cdot(y-x).
\]

To obtain the bounded displacement property, notice that for any $g\in
\supp\mu$
and for any $x\in\R$, the domain that we have used to define $\Phi
_g(x)$ contains
[as $g$ is $1/\mu(g)$-Lipschitz] a rectangular triangle with sides
$|x-g(x)|$ and
$\mu(g_0)\cdot|x-g(x)|$; see Figure~\ref{figTriangles}.

%f2 #&#
%
\begin{figure}

\includegraphics{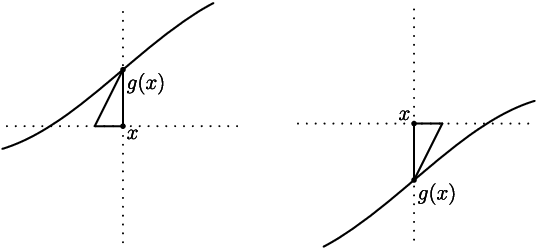}

\caption{Triangles.}\label{figTriangles}
\end{figure}

Hence, $\Phi_g(x)\ge\frac{\mu(g)}{2} |x-g(x)|^2$, which implies that
$\Phi_\mu\ge\frac{\mu(g)^2}{2} |x-g(x)|^2$. Since $\Phi_\mu$ does
not depend on $x$, we obtain the desired uniform upper bound for
the displacement $|g(x)-x|$.
\end{pf}

As a consequence, we have the following result.
%
%th8.5 #&#
%
\begin{theorem}\label{lip}
If $G$ is an irreducible, finitely generated subgroup of $\mathrm
{Homeo}^+({\mathbf R})$,
then there exists a homeomorphism $D\dvtx {\mathbf R} \rightarrow
{\mathbf R}$ such
that, for
every $g\in G$, the map $D \circ g \circ D^{-1}$ is Lipschitz and has uniformly
bounded displacement.
\end{theorem}
\begin{pf} Without loss of generality, assume that the $G$-action is
minimal---otherwise,
consider the subgroup of $\mathrm{Homeo}^+({\mathbf R})$ generated by $G$
and two rationally
independent translations. By Theorem~\ref{Tderrienniccoordinates},
$G$ is semi-conjugate
to a group satisfying the Derriennic property. Since the orbits of $G$
are dense, the
semiconjugacy is in fact a conjugacy. Now the desired conclusion
follows as an
application of Proposition~\ref{Llipschitz}.\vadjust{\goodbreak}
\end{pf}

It should be pointed out that the result above also follows from
\cite{DKN}, Theorem~D
(see also~\cite{Deroin}) by means of rather tricky---and less
conceptual---arguments.
The reader is referred to~\cite{book} for a detailed discussion on
this. Finally, a
conjugacy into a group of $C^1$ diffeomorphisms of the line is not
always possible;
see~\cite{GT} and references therein.

\section*{Acknowledgments}

The authors would like to thank A. Erschler and S. Brofferio for many
interesting and stimulating conversations.
K. Parwani would like to thank the people at the D\'epartement de
Math\'ematiques d'Orsay for their gracious hospitality during the visit
in which this project was initiated.

%suskaldyti doi

% imsref loaded by lrinkeviciute, 2013-01-03 09:28:26
%

\printaddresses

\end{document}